
\input amstex
\loadeusm
\documentstyle{amsppt}
\magnification=\magstep{1}  
\pagewidth{6.5truein}
\pageheight{9.0truein}
\NoBlackBoxes

\long\def\ignore#1\endignore{#1}

\input xy \xyoption{matrix} \xyoption{arrow}
          \xyoption{curve}  \xyoption{frame}
\def\edge{\ar@{-}}
\def\eddge{\ar@{=}}
\def\dttdar{\ar@{.>}}
\def\drbl{\save+<0ex,-2ex> \drop{\bullet} \restore}

\def\dashedge{\ar@{--}}
\def\ddashedge{\ar@{==}}
\def\dotedge{\ar@{..}}

\def\dropleft#1#2{\save+<-#1ex,0ex>\drop{#2}\restore}
\def\dropup#1#2{\save+<0ex,#1ex> \drop{#2} \restore}

\def\la{\Lambda}
\def\lamod{\Lambda{}\operatorname{-mod}}

\def\mod{\operatorname{mod}}
\def\Mod{\operatorname{Mod}}
\def\pinflamod{\operatorname{\P^{< \infty}(\lamod)}}

\def\pinf{\operatorname{\P^{< \infty}}}
\def\pdim{\operatorname{p\,dim}}

\def\Aut{\operatorname{Aut}}
\def\findim{\operatorname{fin\,dim}}
\def\Findim{\operatorname{Fin\,dim}}
\def\repdim{\operatorname{rep\,dim}}

\def\soc{\operatorname{soc}}
\def\supp{\operatorname{supp}}
\def\St{\operatorname{St}}
\def\Bd{\operatorname{Bd}}

\def\NN{{\Bbb N}}
\def\ZZ{{\Bbb Z}}

\def\P{{\Cal P}}

\def\la{\Lambda}
\def\lamod{\Lambda{}\operatorname{-mod}}

\def\pinflamod{\operatorname{\P^{< \infty}(\lamod)}}

\def\pdim{\operatorname{p\,dim}}

\def\Aut{\operatorname{Aut}}
\def\findim{\operatorname{fin\,dim}}
\def\Findim{\operatorname{Fin\,dim}}

\def\soc{\operatorname{soc}}

\def\NN{{\Bbb N}}
\def\ZZ{{\Bbb Z}}

\def\P{{\Cal P}}

\def\BHT{{\bf  1}}
\def\Bass{{\bf 2}}
\def\Bon{{\bf 3}}
\def\BuRi{{\bf 4}}
\def\DoFr{{\bf 5}}
\def\Erd{{\bf 6}}
\def\EHIS{{\bf 7}}
\def\Gal{{\bf 8}}
\def\GePo{{\bf 9}}
\def\pre{{\bf 10}}
\def\dom{{\bf 11}}
\def\SmaHZ{{\bf 12}}
\def\SmaHZstring{{\bf 13}}
\def\IgTo{{\bf 14}}
\def\Kra{{\bf 15}}
\def\Rin{{\bf 16}}
\def\Smalo{{\bf 17}}
\def\WaWa{{\bf 18}}

\topmatter

\title Representation-tame algebras need not be homologically tame \endtitle

\author Birge Huisgen-Zimmermann \endauthor

\abstract  We show that, also within the class of representation-tame finite dimensional algebras $\la$, the big left finitistic dimension of $\la$ may be strictly larger than the little.  In fact, the discrepancies $\, \Findim \la - \findim \la$ need not even be bounded for special biserial algebras which constitute one of the (otherwise) most thoroughly understood classes of tame algebras.  

More precisely:  For every positive integer $r$, we construct a special biserial algebra $\la$ with the property that $\findim \la = r + 1$, while $\Findim \la = 2r + 1$.  In particular, there are infinite dimensional representations of $\la$ which have finite projective dimension, while not being direct limits of {\it finitely generated\/} representations of finite projective dimension.

\endabstract

\thanks This work was partially supported by a grant from the National Science Foundation. \endthanks

\address Department of Mathematics, University of California, Santa
Barbara, CA 93106-3080 \endaddress
\email birge\@math.ucsb.edu \endemail

\endtopmatter

\document

\head 1. Introduction and conventions \endhead

Special biserial algebras have a particularly transparent tame finite dimensional representation theory.  This was shown by Gelfand and Ponomarev for a subclass which is representation-theoretically linked to the Lorentz group \cite{\GePo}.  Their classification of the finite dimensional representations in a restricted scenario was  incrementally extended by several authors, the final step being due to 
Wald and Waschb\"usch \cite{\WaWa} (see Theorem 0 below).  Our primary goal here is to show that the infinite dimensional representation theory of special biserial algebras may, from a homological viewpoint, ``wildly" deviate from the finite dimensional:  Namely, for any positive integer $m$, there is a special biserial algebra $\la$ with the property that $\Findim \la - \findim \la \ge m$; here $\findim \la$ and $\Findim \la$ stand for the left little and big finitistic dimensions of $\la$.  More precisely, for each $r \ge 1$, there is a special biserial algebra $\la$ with the property that $\findim \la = r + 1$, while $\Findim \la = 2r +1$.  In particular, the structure of infinite dimensional $\la$-modules of finite projective dimension is not simply based on amalgamations of finite dimensional patterns.

On the other hand, we conjecture that, for special biserial algebras $\la$, the big finitistic dimension cannot grow any faster as a function of the little than in our examples;  in other words, we conjecture that $\, \Findim \la \le  2 \findim \la - 1$.  If confirmed, the conjecture clearly implies finiteness of $\Findim \la$, given that $\findim \la$ is already known to be finite.
This latter fact was proved by Erdmann, Holm, Iyama, and Schr\"oer in \cite{\EHIS}, where it arose as a consequence of the result that the representation dimension of a special biserial algebra is at most $3$; indeed, by work of Igusa and Todorov \cite{\IgTo}, the inequality $\repdim \la \le 3$ implies finiteness of the little finitistic dimension. However, the big finitistic dimension eludes this approach.  

To fill in some background:  The first finitistic dimension conjecture, originally stated as a problem in \cite{\Bass}, postulated equality of $\findim \la$ and $\Findim \la$ whenever $\la$ is finite dimensional.  Its failure first surfaced in \cite{\dom}, where it was shown that the two dimensions may differ even for monomial algebras.  While the discrepancy $\, \Findim \la - \findim \la$ cannot exceed $1$ in the monomial case, it is known to take arbitrarily large values for more general classes of wild algebras.  Illustrations of this phenomenon are obtained from the monomial examples in \cite{\dom} by way of Rickard's observation that both big and little finitistic dimensions behave additively on tensor products $\la_1 \otimes_K \la_2$; a more direct construction was presented by Smal\o \ in \cite{\Smalo}. 
On the positive side, the equality $\Findim \la  = \findim \la$ has been secured for large classes of algebras.  For example, it holds whenever the category $\pinflamod$ is contravariantly finite in $\lamod$ (\cite{\SmaHZ}), this being just the tip of the iceberg of positive instances.   

That the homological behavior of special biserial algebras should be understood so late in the game, and in slow increments at that, is somewhat surprising, as they constitute one of the most thoroughly investigated classes of tame algebras, next to the hereditary algebras based on (extended) Dynkin graphs.  They have, in fact, developed into a  showcase for representation-theoretic techniques, due to the combined facts that $\bullet$ they occur widely in contexts of interest (such as the representation theory of the Lorentz group and among blocks of group algebras in characteristic $2$;  see \cite{\GePo}, \cite{\Bon}, \cite{\Rin}, \cite{\DoFr} and \cite{\Erd}, for instance) and $\bullet$ the structure of their indecomposable finite dimensional representations is fully understood; it is governed by two simple templates, strings and bands (for use in our arguments, we define them below).  The firm grip on the finite dimensional representations was, in turn, extensively used towards understanding Auslander-Reiten quivers and numerous other aspects of special biserial algebras, while the homological analysis lagged behind.   The first foray in the latter direction, in \cite{\SmaHZstring}, addressed only string algebras, i.e., special biserial algebras which are monomial; in this setting, the two finitistic dimensions do coincide, and the finite value is easily computed from quiver and relations.   It was followed by the seminal result in \cite{\EHIS} regarding the little finitistic dimension in the general special biserial case.

\bigskip
{\it Notation and Terminology.} Our (arbitrary) base field will be consistently labeled $K$.  For any finite dimensional algebra $\Delta$, we denote by $\Delta$-$\Mod$ and $\Delta$-$\mod$ the category of all left $\Delta$-modules and the subcategory consisting of the finitely generated left $\Delta$-modules, respectively.  By $\pinf(\Delta$-$\mod)$, resp\. $\pinf(\Delta$-$\Mod)$, we mean the full subcategory based on the objects of finite projective dimension in the indicated module category.   The {\it left little and big finitistic dimensions\/} of $\Delta$ are
$$\findim \Delta \ \ =\ \  \sup\{\pdim M \mid M \in \pinf(\Delta\text{-}\mod)\}$$
and 
$$\Findim \Delta \ \ =\ \  \sup\{\pdim M \mid M \in \pinf(\Delta\text{-}\Mod)\}.$$
Even though both of these invariants are side-sensitive, we will suppress the qualifier ``left", since we will focus on left modules throughout.  

All of our algebras $\Delta$ will be path algebras modulo admissible ideals.  If $\Delta = KQ/I$, we identify the set of vertices of $Q$ with a full set $\{e_1, \dots, e_n\}$ of orthogonal primitive idempotents of $\Delta$.  Moreover $J = J(\Delta)$ will denote the Jacobson radical. Given $M \in \Delta$-$\Mod$, we call an element $x \in M$ a {\it top element\/} of $M$ if $x \in M \setminus JM$ and $x$ is normed by some vertex $e_i$, i.e., $x = e_i x$.  A {\it full family of top elements\/} of $M$ is a family $(x_a)_{a \in A}$ of top elements which generates $M$ and is $K$-linearly independent modulo $JM$.  Note:  Given such a full family of top elements of $M$, the direct sum $\bigoplus_{1 \le i \le n} (\la e_i)^{(A_i)}$ is a projective cover of  both $M$ and $M/JM$; here $A_i = \{a \in A \mid e_i x_a = x_a\}$.   Our convention for composing paths in $KQ$ is as follows: $qp$ stands for ``q after p".
\medskip

\noindent {\it Background on special biserial algebras.\/}

Throughout, $\la$ will stand for a {\it special biserial algebra\/}.  This means that $\la = KQ/I$, where $\bullet$ $Q$ is a quiver with the property that no vertex occurs as the starting point of more than two arrows or as the endpoint of more than two arrows; and $\bullet$ for every arrow $\alpha$ of $Q$, there is at most one arrow $\beta$ with $\alpha \beta \notin I$ and at most one arrow $\gamma$ with $\gamma \alpha \notin I$.  We refer to the extensive bibliography of \cite{\SmaHZstring} for much of the work on this class of algebras.

The definition impinges on the structure of the indecomposable projective left $\la$-modules $\la e_i$ as follows:  Either {\bf(a)} $J e_i$ is the direct sum of two uniserial modules $U_1$ and $U_2$ (with $U_j = 0$ permissible), or else {\bf(b)} $J e_i = U_1 + U_2$ with $U_j$ uniserial and $U_1 \cap U_2$ simple. 

We start by recalling the pivotal theorem that classifies the indecomposable objects in $\lamod$.
 In its present form, it was established by Wald and Waschb\"usch \cite{\WaWa}, the underlying ideas having evolved in a sequence of successive generalizations.

\proclaim{Theorem 0} {\rm (See \cite{\GePo, \Rin, \Bon, \DoFr,
\BuRi, \WaWa})} Apart from the indecomposable projective left $\la$-modules of type {\rm{(}}b{\rm{) }}above, the indecomposable representations in $\lamod$ are either string or band modules {\rm{(}}to be described next{\rm{)}}.  Conversely, all strings and bands are indecomposable.\qed\endproclaim

We slightly modify the existing notation to describe string and band
modules for our present purpose.  In particular, we address their graphs (in the sense of \cite{\BHT, Definition 3.9}), since those provide the most convenient means of computing and displaying syzygies. The set $\Cal P$ of paths that feed into the definition of ``words" depends on both $Q$ and $I$:  We call a path $p$ in $KQ \setminus I$ a {\it syllable\/} in case $p \ne kq$ modulo $I$ for all $k \in K^*$ and all paths $q \in KQ \setminus\{p\}$;  if $p$ starts in the vertex $e_i$, the latter means that the indecomposable projective module $\la e_i$ is either of type (a) above, or else satisfies $\la p \supsetneqq \soc \la e_i$ in type (b).   More generally, a {\it syllable\/} is any element of the set $\Cal P \sqcup {\Cal P}^{-1}$.
  The paths of length $0$, i.e., the vertices of $Q$ will be called the {\it trivial paths}; the trivial
paths and their inverses are also referred to as {\it trivial syllables}.  {\it {\rm{(}}Generalized{\rm{)}} words\/} are
$\ZZ$-indexed sequences of pairs of syllables $w = (p_i^{-1} q_i)_{i \in
\ZZ}$ with $p_i, q_i \in
\Cal P$, which we also communicate as juxtapositions
$$\dots (p_r^{-1} q_r) \dots (p_{-1}^{-1} q_{-1}) (p_0^{-1} q_0) (p_1^{-1} q_1) \dots
(p_s^{-1} q_s) \dots$$ 
 subject to the following constraints:
\roster
\item"$\bullet$"  For each $i \in \ZZ$, the starting points of $p_i$ and
$q_i$ coincide, but the first arrows of $p_i$ and $q_i$ are distinct
whenever both
$p_i$ and $q_i$ are nontrivial.
\item"$\bullet$"  For each $i \in \ZZ$, the end points of $q_i$ and
$p_{i+1}$ coincide, but the last arrows of $q_i$ and $p_{i+1}$ are distinct
whenever both $q_i$ and $p_{i+1}$ are nontrivial.
\item"$\bullet$" No trivial syllables occur between two nontrivial syllables
(i\.e., the nontrivial syllables form a `connected component').
\endroster

\noindent A word  $w = (p_i^{-1} q_i)_{i \in \ZZ}$ will be called {\it finite}
in case, for all $i \gg 0$ and all $i \ll 0$, the syllables with index $i$ are trivial; finite words are also communicated as finite
juxtapositions
$(p_i^{-1} q_i)$ in which the nontrivial syllables are preserved.  More
generally, we do not insist on recording trivial syllables; keep in mind
that they can only occur at the left or right tail ends of a word.  It is
self-explanatory what we mean by a {\it left\/} or {\it right finite\/}
word, and by a {\it left\/} or {\it right\/} periodic generalized word.  

{\it String modules\/}:  Each (generalized) word $w = (p_i^{-1} q_i)_{i \in \ZZ}$ gives rise to a
{\it {\rm{(}}generalized{\rm{)}} string module\/} $\St(w)$, that is, a module $M$ having a graph of the form
$$\xymatrixrowsep{2pc}\xymatrixcolsep{1pc}
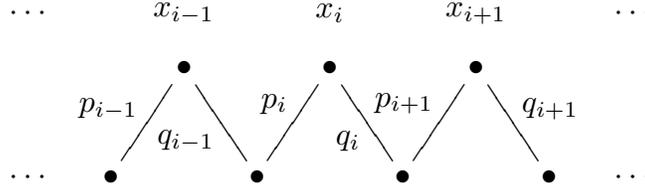
\xymatrix{
\dropup{4}{\cdots} &&\bullet \dropup{4}{x_{i-1}} \edge[dl]_{p_{i-1}}
\edge[dr]_(0.55){q_{i-1}} &&\bullet \dropup{4}{x_i}
\edge[dl]_(0.45){p_i} \edge[dr]_(0.55){q_i} &&\bullet \dropup{4}{x_{i+1}}
\edge[dl]_(0.45){p_{i+1}} \edge[dr]^{q_{i+1}} &&\dropup{4}{\cdots} \\
\cdots &\bullet &&\bullet &&\bullet &&\bullet &\cdots }$$
relative to a full family $(x_i)_{i \in \ZZ}$ of top elements of $M$. 
More formally:  If $w$ is
trivial, say $w = e$, then $\St(w)$ is the simple module $\la e/ Je$.  Now
suppose that $w$ is nontrivial, and let $\supp(w)$ be the set of all those
integers $j$ for which either $p_j$ or $q_j$ is nontrivial. Moreover, let $e(i)$ be
the joint starting vertex of $p_i$ and $q_i$ and $\la z_i$ a copy of $\la e(i)$ with $e(i)z_i = z_i$. Then
$$\St(w) = \biggl(\, \bigoplus_{i \in \supp(w)}\la z_i \biggr)
\biggm/ C \ , \quad \quad \text{where}$$
$$C\ \  = \ \biggl(\, \sum_{i, i+1 \in
\supp(w)} \la \bigl( q_i z_i - p_{i+1} z_{i+1} \bigr)\,\biggr)\ +\
C_{\text{left}} \ + \ C_{\text{right}},$$ 
with cyclic correction terms
$C_{\text{left}}$ and $C_{\text{right}}$ trimming the left and right ends in case $w$ is finite; they are defined as follows:
$C_{\text{left}} = 0$ if either $\supp(w)$ is unbounded on the negative
$\ZZ$-axis or else $l = \inf \supp(w) \in \ZZ$ and there is no arrow $\alpha$ with the property that $\alpha p_l$ is again a path in $KQ \setminus I$; in the remaining case, where $l \in \ZZ$ and there exists an
arrow $\alpha$ (necessarily unique) such that $\alpha
p_l \in KQ \setminus I$, we set
$C_{\text{left}} = \la\, \alpha p_l \, z_l$.  The right-hand correction term
$C_{\text{right}}$ is defined symmetrically.  Clearly, $\St(w)$ is finite
dimensional over $K$ precisely when $w$ is a finite word; the finite dimensional string modules are the ``traditional" ones. That also infinite dimensional (i.e., generalized) string modules are indecomposable was shown in
\cite{\Kra}.  Note moreover that all of the indecomposable projective $\la$-modules which are of type (a), as introduced above, are among the string modules. 
\smallskip

{\it Band modules\/}: The second class of indecomposable representations of $\lamod$ is defined as follows:  Suppose that
$v= p_0^{-1}q_0
\dots p_t^{-1}q_t$ is a finite word with $t\ge 0$ and $p_0$, $q_t$ both
nontrivial; by our conventions, this amounts to the same as to require that
all $p_i$ and $q_i$ be nontrivial. We call
$v$ {\it primitive\/} if
\roster
\item"$\bullet$" the juxtaposition $v^2 =vv$ is again a word (in which case
all powers
$v^r$ are words), and
\item"$\bullet$" $v$ is not itself a power of a strictly shorter word.
\endroster

\noindent In addition to the primitive word $v$,  let $r$ be a positive integer and
$\phi: K^r \rightarrow K^r$ an irreducible automorphism with Frobenius companion matrix
$$\pmatrix 0&\cdots&0&c_1\\ 1&\ddots &&\vdots\\ &\ddots&0&\vdots\\
0&\cdots&1&c_r \endpmatrix.$$  
Moreover, let $\St(v^r)$ be the string module with graph

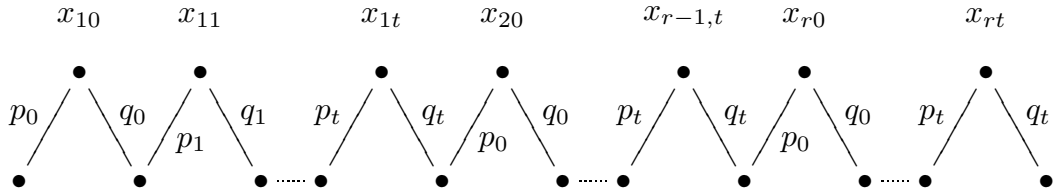
\ignore
$$\xymatrixrowsep{2pc}\xymatrixcolsep{0.67pc}
\xymatrix{
 &\bullet \dropup{4}{x_{10}} \edge[dl]_{p_0} \edge[dr]^{q_0} &&\bullet
\dropup{4}{x_{11}} \edge[dl]^(0.55){p_1}
\edge[dr]^{q_1} &&&\bullet \dropup{4}{x_{1t}}
\edge[dl]_{p_t} \edge[dr]^{q_t} &&\bullet \dropup{4}{x_{20}}
\edge[dl]^(0.55){p_0}
\edge[dr]^{q_0} &&&\bullet \dropup{4}{x_{r-1,t}}
\edge[dl]_{p_t} \edge[dr]^{q_t} &&\bullet \dropup{4}{x_{r0}}
\edge[dl]^(0.55){p_0}
\edge[dr]^{q_0} &&&\bullet  \dropup{4}{x_{rt}}
\edge[dl]_{p_t} \edge[dr]^{q_t}\\
\bullet &&\bullet &&\bullet  \ar@{.}[r] &\bullet &&\bullet &&\bullet
\ar@{.}[r] &\bullet
 &&\bullet &&\bullet  \ar@{.}[r] &\bullet &&\bullet }$$
\endignore

\noindent relative to a suitable full family $x_{10},\dots,x_{1t}$,
$x_{20}, \dots, x_{2t}$, $\dots, x_{r0}, \dots, x_{rt}$ of top elements.  In particular,  $q_jx_{ij}= p_{j+1}x_{i,j+1}$ for $0\le i\le r$ and
$1\le j<t$, and
$q_tx_{it}= p_0x_{i+1,0}$ for $i<r$.  Then the {\it band module\/}
$\Bd(v^r,\phi)$ is defined as follows: 
$$\Bd(v^r,\phi)= \St(v^r)/ \la \bigl( q_t x_{rt}- \sum_{i=1}^r c_ip_0x_{i0}
\bigr).$$ 
Clearly, the canonical images $y_{ij}$ of the $x_{ij}$ in
$\Bd(v^r,\phi)$ constitute a full family of top elements of
$\Bd(v^r,\phi)$; they
satisfy the string equations, as well as the one additional equation
$$q_ty_{rt}= \sum_{i=1}^r c_i p_0y_{i0}.$$  
Relative to the family $(y_{ij})$, the band module
$\Bd(v^r,\phi)$ is graphed in the form

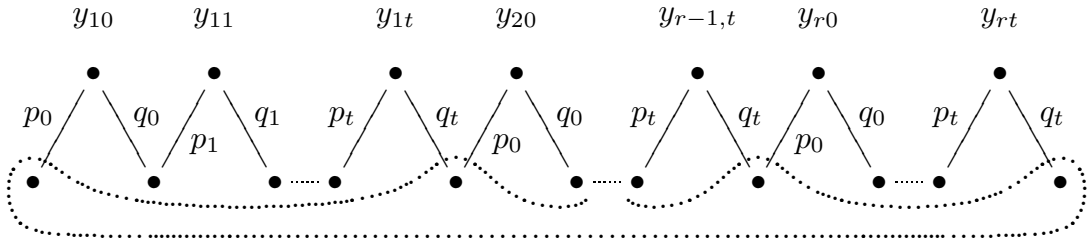
\ignore
$$\xymatrixrowsep{2pc}\xymatrixcolsep{0.67pc}
\xymatrix{
 &\bullet \dropup{4}{y_{10}} \edge[dl]_{p_0} \edge[dr]^{q_0} &&\bullet
\dropup{4}{y_{11}} \edge[dl]^(0.55){p_1}
\edge[dr]^{q_1} &&&\bullet \dropup{4}{y_{1t}}
\edge[dl]_{p_t} \edge[dr]^{q_t} &&\bullet \dropup{4}{y_{20}}
\edge[dl]^(0.55){p_0}
\edge[dr]^{q_0} &&&\bullet \dropup{4}{y_{r-1,t}}
\edge[dl]_{p_t} \edge[dr]^{q_t} &&\bullet \dropup{4}{y_{r0}}
\edge[dl]^(0.55){p_0}
\edge[dr]^{q_0} &&&\bullet  \dropup{4}{y_{rt}}
\edge[dl]_{p_t} \edge[dr]^{q_t}\\
\bullet &&\bullet &&\bullet  \ar@{.}[r] &\bullet &&\bullet &&\bullet
\ar@{.}[r] &\bullet
\save[0,0]+(-1,-2);[0,-1]+(1,-2) **\crv{~*=<2.5pt>{.} [0,0]+(0,-3)
&[0,1]+(0,-3) &[0,2]+(-2,3) &[0,2]+(2,3) &[0,3]+(0,-3) &[0,6]+(0,-3)
&[0,7]+(-2,3) &[0,7]+(3,3) &[0,7]+(3,-6) &[0,6]+(0,-6) &[0,-9]+(0,-6)
&[0,-10]+(-3,-6) &[0,-10]+(-3,3) &[0,-10]+(2,3) &[0,-9]+(0,-3)
&[0,-4]+(0,-3) &[0,-3]+(-2,3) &[0,-3]+(2,3) &[0,-2]+(0,-3) &[0,-1]+(0,-3)
}\restore
 &&\bullet &&\bullet  \ar@{.}[r] &\bullet &&\bullet }$$
\endignore
\vskip0.3truein

\noindent where the dotted line that encircles the vertices representing the
elements $p_0 y_{10},\dots,
\allowmathbreak p_0y_{r0}, q_ty_{rt}$ in the socle of $\Bd(v^r,\phi)$ encodes the information that 
the $K$-space spanned by these $r+1$ elements has dimension $r$.

In contrast to the graph of a string module, the graph of a
band module $B = \Bd(v^r,\phi)$ does not pin $B$ down up to isomorphism, but  communicates a one-parameter family of band modules.  It is the combination of the graph with the eigenvalue of $\phi$ that determines the isomorphism class of a specific member of this family.  

\medskip

\head 2.  Syzygies over special biserial algebras  \endhead

As before, $\la = KQ/I$ is a special biserial algebra.  The folllowing auxiliary facts concerning syzygies and projective dimensions of band and (generalized) string modules over special biserial algebras will be freely used in the next section.  
The arguments backing them can be found in Julia Galstad's thesis (\cite{\Gal}), which is presently in progress.  

\proclaim{Proposition 2.1}  Any syzygy of a {\rm{(}}generalized{\rm{)}} string module is a direct sum of generalized string modules. Moreover, if $w = (p_i^{-1}q_i)_{i \in \ZZ}$ is a generalized word, then the following are equivalent:
\smallskip
{\bf (a)} $\Omega^1\bigl(\St(w)\bigr)$ is decomposable.

 {\bf (a)} The projective cover of $\St(w)$ contains at least one direct summand that is a 
string 

\quad \ \  module {\rm{(}}in other words, at least one of the indecomposable direct summands of this 

\quad \ \ cover is projective of type {\rm{(}}a{\rm{)}}{\rm{)}}.
 \endproclaim 

\proclaim{Proposition 2.2}   Let $v$ be a primitive word and  $\widehat{v}$ the corresponding two-sided infinite generalized word $\, ... vvv...\, $.  Moreover, let $m$ be any positive integer and $\psi$ an irreducible automorphism of $K^m$.  Then:
\smallskip

{\bf (a)}  The following statements are equivalent:

$\bullet$ $\Omega^1\bigl(\Bd(v^m, \psi) \bigr)$ is a band module.

$\bullet$ $\Omega^1\bigl(\St(\widehat{v})\bigr)$ is an indecomposable generalized string module. 

$\bullet$  None of the indecomposable direct summands of the projective cover of $\St(v)$ is 

\quad a string module.
\smallskip

If these conditions fail to be satisfied, then $\Omega^1\bigl(\Bd(v^m, \psi) \bigr)$ is a direct sum of string modules.
\smallskip
 
 {\bf (b)} For any $m \ge 1$ and any cyclic automorphism $\psi \in \Aut_K (K^m)$, 
 $$\pdim \Bd(v^m, \psi) = \pdim  \St(\widehat{v}).$$
 \smallskip
 
 {\bf (c)} If $u v v v \dots$ is a right periodic generalized word such that $\St(u v v v \dots)$ has finite projective dimension, then $\pdim \Bd(v^m , \psi) < \infty$ for any $m \ge 1$ and any cyclic automorphism $\psi$ of $K^m$.  A symmetric statement holds for left periodic generalized words.
 \endproclaim
 
 In the upcoming example, the big finitistic dimension of $\la$ will be attained on a right periodic generalized string module as described in part (c) of the preceding proposition.  
 We note that, for any band module $B$, a graph of $\Omega^1(B)$ is obtainable at a glance from a  graph of $B$ and the graphs of the indecomposable projective modules $\la e_i$.  That the same holds for string modules is still more obvious.  Illustrations will abound in Section 3.

\head  3. The family of examples \endhead

Fix $r \ge 1$.  We begin by constructing a string algebra (a special biserial algebra which is also monomial), labeled $\la_1 = KQ^{(1)}/ I_1$, with $\findim \la_1 = \Findim \la_1 = r + 1$, and then recursively move to successive one-point extensions.  This process will yield special biserial algebras $\la_m$ for all $m$ with the property that $\Findim \la_m = \Findim \la_{m - 1} + 1$, while the little finitistic dimensions of the $\la_m$ remain stationary for $m \le r+1$.   The algebra $\la_{r+1}$ will then realize the maximal discrepancy between big little finitistic dimensions occurring in the series $\la_1, \la_2, \la_3, \dots$; namely $\Findim \la_{r+1} = 2r + 1 = 2 \findim \la_{r+1} - 1$.
 
The quiver $Q^{(1)}$ is as follows, and the paths generating $I_1$ can be gleaned from the graphs of the indecomposable projective left $\la_1$-modules given below. 

$$\xymatrixcolsep{1pc}\xymatrixrowsep{2pc}
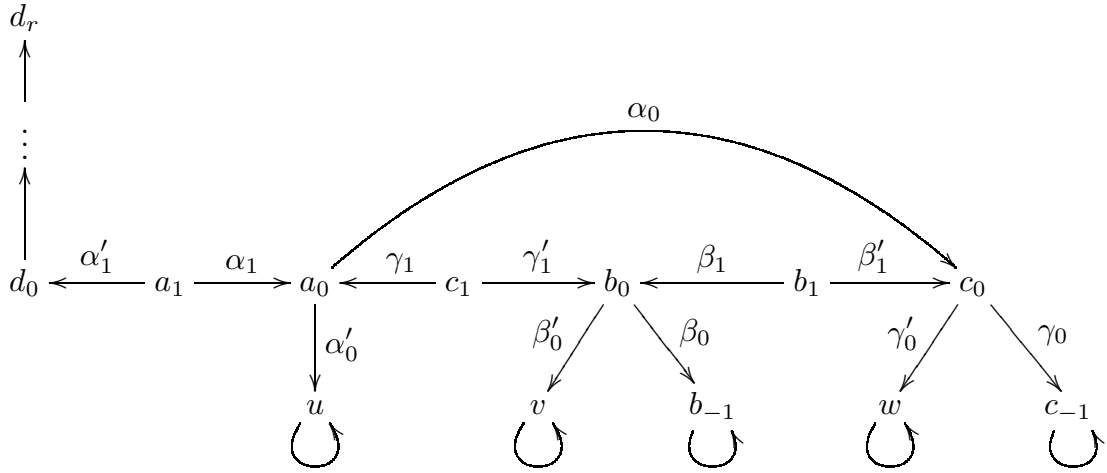
\xymatrix{
d_r  \\
\vdots \ar[u]  \\
d_0 \ar[u] &&a_1 \ar[ll]_{\alpha'_1} \ar[rr]^{\alpha_1} &&a_0 \ar@/^4pc/[0,8]^{\alpha_0} \ar[d]^{\alpha'_0} &&c_1 \ar[ll]_(0.4){\gamma_1} \ar[rr]^{\gamma'_1} &&b_0 \ar[dl]_(0.55){\beta'_0} \ar[dr]^(0.55){\beta_0} &&b_1 \ar[ll]_{\beta_1} \ar[rr]^(0.4){\beta'_1} &&c_0 \ar[dl]_(0.55){\gamma'_0} \ar[dr]^(0.55){\gamma_0}  \\
 &&&&u \ar@(dl,dr) &&&v \ar@(dl,dr) &&b_{-1} \ar@(dl,dr) &&w \ar@(dl,dr) &&c_{-1} \ar@(dl,dr) \
 }$$
 \bigskip
 \bigskip
\noindent  The graphs of the indecomposable projectives in $\la_1$-mod:

$$\xymatrixcolsep{0.1pc}\xymatrixrowsep{1.5pc}
\xymatrix{
 &a_0 \edge[dl] \edge[dr] &&&&a_1 \edge[dl] \edge[dr] &&&&b_0 \edge[dl] \edge[dr] &&&&b_1 \edge[dl] \edge[dr] &&&&c_0 \edge[dl] \edge[dr] &&&&c_1 \edge[dl] \edge[dr]  \\
u &&c_0 \edge[d] &&d_0 &&a_0 \edge[d] &&v &&b_{-1} &&b_0 \edge[d] &&c_0 \edge[d] &&w &&c_{-1} &&a_0 \edge[d] &&b_0 \edge[d]  \\
 &&c_{-1} &&&&u &&&&&&b_{-1} &&w &&&&&&c_0 \edge[d] &&v  \\
 &&&&&& &&&&&&&& &&&&&&c_{-1}  \\
b_{-1} \edge[d] &&c_{-1} \edge[d] &&u \edge[d] &&v \edge[d] &&w \edge[d] &&d_0 \edge[d] &&d_1 \edge[d] &\cdots &d_{r-1} \edge[d] &&d_r \drbl  \\
b_{-1} &&c_{-1} &&u &&v &&w &&d_1 &&d_2 &\cdots &d_r
}$$
\bigskip

\noindent {\bf Claim 1.}  $\findim \la_1 = \Findim \la_1 = r+1$,   and 
$$r =  \max\{ \la_1\, q \mid q \ \text{is a path of positive length in}\  KQ^{(1)}\setminus I_1 \ \text{such that}\  \pdim \la_1\, q < \infty\}.$$
Moreover, the only path $p$ with $\pdim \la_1 p = r$ is $p = \alpha'_1$.  

\demo{Proof}  That $r$ equals the displayed maximum is straightforward.  Using \cite{\pre, Corollary II}, we deduce that $\findim \la_1$ and $\Findim \la_1$ both belong to $\{ r+1, r+2\}$.  Given that $\alpha'_1$ is the only path generating a cyclic module of projective dimension $r$, we derive from  \cite{\pre, Theorem VI} that, in fact, $\findim \la = \Findim \la = r+1$.   
\qed \enddemo 

The construction of $\la_m$ for $m =2$ deviates from the pattern that underlies the definition of the $\la_m$ for $m \ge 3$ and will be described separately.  The quiver $Q^{(2)}$ is obtained from $Q^{(1)}$ through three successive one-point extensions, the first adding a vertex $c_2$, followed by the addition of two further vertices $a_2$, $b_2$; there are six new arrows as shown below. 

$$\xymatrixcolsep{0.5pc}\xymatrixrowsep{1.75pc}
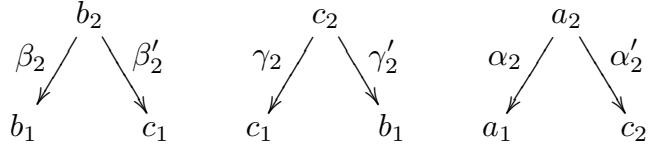
\xymatrix{
 &b_2 \ar[dl]_{\beta_2} \ar[dr]^{\beta'_2} &&&&c_2 \ar[dl]_{\gamma_2} \ar[dr]^{\gamma'_2} &&&&a_2 \ar[dl]_{\alpha_2} \ar[dr]^{\alpha'_2}  \\
b_1 &&c_1 &&c_1 &&b_1 &&a_1 &&c_2
}$$

\noindent We define $\la_2 = KQ^{(2)} / I_2$, where $I_2$ is the ideal generated by $I_1$ and the following elements  in $KQ^{(2)}$: $ \alpha_1 \alpha_2 - \gamma_1 \gamma_2 \alpha'_2$, $\gamma'_1 \beta'_2 - \beta_1 \beta_2$, and $\beta'_1 \gamma'_2 - \alpha_0 \gamma_1 \gamma_2$, next to all those paths starting in one of the vertices $a_2$, $b_2$, $c_2$ which are not subpaths of any of the six paths involved in the preceding binomial relations.  The relations imposed on $\la_2$ are more suggestively communicated (if only up to nonzero scalars) by the graphs of the three additional indecomposable projective left $\la_2$-modules: 

$$\xymatrixcolsep{0.75pc}\xymatrixrowsep{1.5pc}
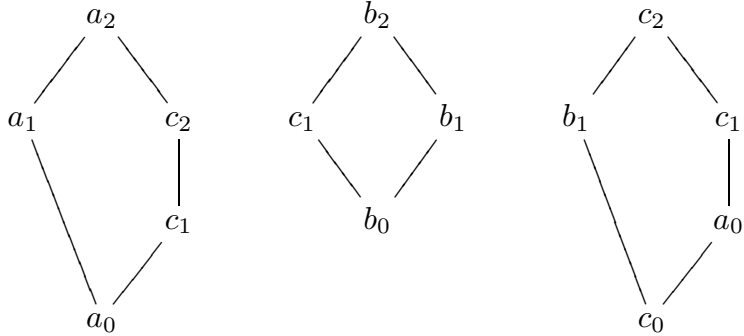
\xymatrix{
 &a_2 \edge[dl] \edge[dr] &&&&b_2 \edge[dl] \edge[dr] &&&&c_2 \edge[dl] \edge[dr]  \\
a_1 \edge[ddr] &&c_2 \edge[d] &&c_1 \edge[dr] &&b_1 \edge[dl] &&b_1 \edge[ddr] &&c_1 \edge[d]  \\
 &&c_1 \edge[dl] &&&b_0 &&&&&a_0 \edge[dl]  \\
 &a_0 &&&& &&&&c_0
}$$

\medskip
Given $\la_{m-1} = KQ^{(m-1)}/I_{m-1}$ for $m \ge 3$, we define $Q^{(m)}$ to be the quiver obtained from $Q^{(m-1)}$ via two one-point extensions (not interlinked), adding vertices $a_m$, $b_m$, together with four arrows as displayed below.
  
$$\xymatrixcolsep{0.5pc}\xymatrixrowsep{1.75pc}
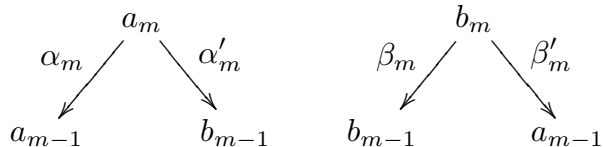
\xymatrix{
 &a_m \ar[dl]_{\alpha_m} \ar[dr]^{\alpha'_m} &&&&b_m \ar[dl]_{\beta_m} \ar[dr]^{\beta'_m}  \\
a_{m-1} &&b_{m-1} &&b_{m-1} &&a_{m-1}
}$$

\noindent  For each integer $m \ge 3$, the ideal $I_m$ is to be generated by $I_{m-1}$, two additional binomial relations, to be described separately for $m \in \{3,4\}$ and $m \ge 5$, and all those paths starting in $a_m$ or $b_m$ which are not subpaths of any of the four paths involved in the added binomial relations in $I_m$.  For $I_3$, the additional binomial relations are 
$\gamma_2 \alpha'_2 \alpha_3 - \beta'_2 \alpha'_3$ and $\beta_2 \beta_3 - \gamma'_2 \beta'_3$.  For $I_4$, they are $\alpha'_3 \alpha_4 - \beta_3 \alpha'_4$ and $\beta'_3 \beta_4 - \alpha'_2 \beta_3 \alpha'_4$.   In a visually intuitive format, these relations are reflected by the graphs of the new indecomposable projective modules $\la_3 e_{a_3}$ and $\la_3 e_{b_3}$, resp., $\la_4 e_{a_4}$ and $\la_4 e_{b_4}$:

$$\xymatrixcolsep{0.75pc}\xymatrixrowsep{1.5pc}
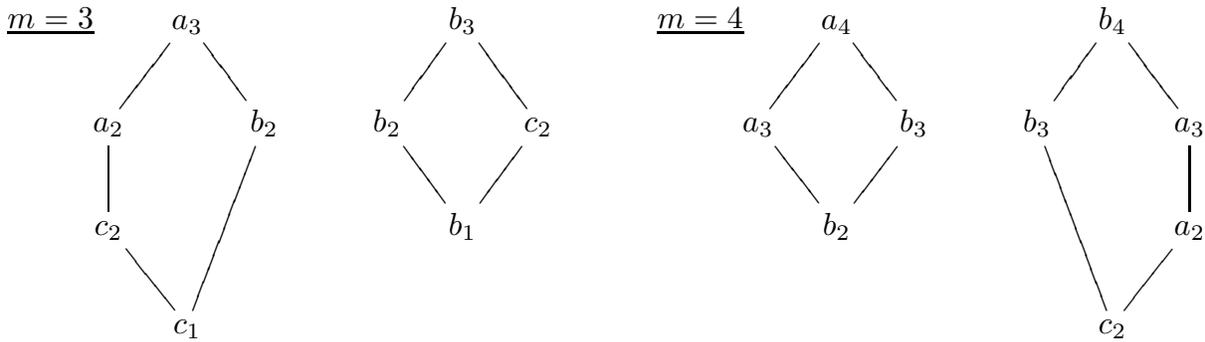
\xymatrix{
 &a_3 \dropleft{10}{\underline{m=3}} \edge[dl] \edge[dr] &&&&b_3 \edge[dl] \edge[dr] &&& &&&a_4 \dropleft{10}{\underline{m=4}} \edge[dl] \edge[dr] &&&&b_4 \edge[dl] \edge[dr]  \\
 a_2 \edge[d] &&b_2 \edge[ddl] &&b_2 \edge[dr] &&c_2 \edge[dl] &&&&a_3 \edge[dr] &&b_3  \edge[dl] &&b_3\edge[ddr] &&a_3 \edge[d]  \\
 c_2 \edge[dr] &&&&&b_1 &&&&&&b_2 &&&&&a_2 \edge[dl]  \\
 &c_1 &&&&& &&&&&&&&&c_2
 }$$   
   
\noindent For $m \ge 5$, finally, the additional relations are $\alpha'_{m-1} \alpha_m -\beta_{m-1} \alpha'_m$ and $\beta'_{m-1} \beta_m - \alpha_{m-1} \beta'_m$.  They yield the following graphs of the indecomposable projective $\la_m$-modules, $\la_m e_{a_m}$ and $\la_m e_{b_m}$:

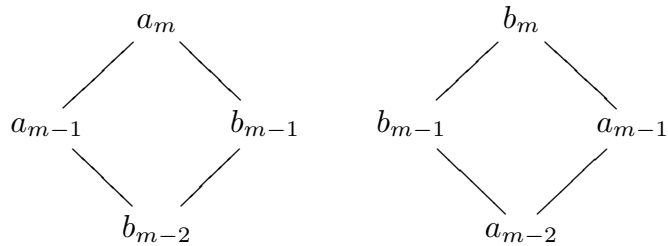
$$\xymatrixcolsep{0.5pc}\xymatrixrowsep{1.5pc}
\xymatrix{
 &a_m \edge[dl] \edge[dr] &&&&b_m \edge[dl] \edge[dr]  \\
a_{m-1} \edge[dr] &&b_{m-1} \edge[dl] &&b_{m-1} \edge[dr] &&a_{m-1} \edge[dl]  \\
 &b_{m-2} &&&&a_{m-2}
}$$

\medskip

Since numerous syzygies need to be computed in ascertaining that this class of algebras behaves as postulated, we include a diagrammatic overview which allows us to find syzygies of (generalized) string modules at a glance; see the {\it Orientation Diagram\/} below.  For any word $w$ ``starting" in a vertex $x_m$, where $x \in \{a,b,c\}$ and $m \ge 2$, either $\St(w)$ or $\St(w^{-1})$ is on display in a horizontal zigzag line, with the corresponding syzygy showing below it (alternately traced with single or double edges).

$$\xymatrixcolsep{1pc}\xymatrixrowsep{3.2pc}
\xymatrix@!0{
 && && && && && && && &&&& &&&& &&&& &&&&  \\
&& && && && && &&a_6 \dotedge[urr] \ddashedge[dll]_{\alpha_6} \eddge[drr] &&&&b_6 \dotedge[urr] \eddge[dll] \eddge[drr] &&&&a_6 \dotedge[urr] \eddge[dll] \eddge[drr] &&&&b_6 \dotedge[urr] \eddge[dll] \eddge[drr] &&&&a_6 \dotedge[urr] \eddge[dll]  \dotedge[drr]  \\
&& && && && &&a_5 \ddashedge[dll]_{\alpha_5} \edge[drr]_(0.6){\alpha'_5} &&&&b_5 \edge[dll]_{\beta_5} \edge[drr]^{\beta'_5} &&&&a_5 \edge[dll] \edge[drr] &&&&b_5 \edge[dll] \edge[drr] &&&&a_5 \edge[dll] \edge[drr] &&&&  \\
&& && && &&a_4 \ddashedge[dll]_{\alpha_4} \eddge[drr]^{\alpha'_4} &&&&b_4 \eddge[dll]^(0.6){\beta_4} \eddge[drr]^{\beta'_4} &&&&a_4 \eddge[dll]^(0.6){\alpha_4} \eddge[drr] &&&&b_4 \eddge[dll] \eddge[drr] &&&&a_4 \eddge[dll] \eddge[drr] &&&&b_4 \dotedge[urr] \eddge[dll] \dotedge[drr] &&&&  \\
&& && &&a_3 \ddashedge[ddll]_{\alpha_3} \edge[ddrr]^(0.4){\alpha'_3} &&&&b_3 \edge[ddll]^(0.6){\beta_3} \edge[ddrr]^(0.4){\beta'_3} &&&&a_3 \edge[dl]_{\alpha_3} \edge[ddrr]^(0.4){\alpha'_3} &&&&b_3 \edge[ddll]^(0.6){\gamma_3} \edge[ddrr]^(0.4){\gamma'_3} &&&&a_3 \edge[dl] \edge[ddrr] &&&&b_3 \edge[ddll] \edge[ddrr] &&&&  \\
 && && && &&&& &&&a_2 \edge[dl]^{\alpha'_2} &&&&& &&&a_2 \edge[dl] &&&&&&&&&  \\
 && &&a_2 \ddashedge[ddll]_{\alpha_2} \eddge[dr]^(0.6){\alpha'_2} &&&&b_2 \eddge[ddll]^{\beta'_2} \eddge[drr]^(0.4){\beta_2} &&&&c_2 \eddge[dll]^(0.6){\gamma'_2} \eddge[drr]^(0.4){\gamma_2} &&&&b_2 \eddge[dll] \eddge[drr] &&&&c_2 \eddge[dll] \eddge[drr] &&&&b_2 \eddge[dll] \eddge[drr] &&&&c_2 \dotedge[urr] \eddge[dll] \dotedge[drr]  \\
 && &&&c_2 \eddge[dr]_{\gamma_2} &&& &&b_1 \edge[ddll]_{\beta_1} \edge[ddrr]_(0.6){\beta'_1} &&&&c_1 \edge[dl]_{\gamma_1} \edge[ddrr]^(0.4){\gamma'_1} &&&&b_1 \edge[ddll] \edge[ddrr] &&&&c_1 \edge[dl] \edge[ddrr] &&&&b_1 \edge[ddll] \edge[ddrr] &&&&  \\
 &&a_1 \ddashedge[dll]_{\alpha'_1} \edge[drr]_(0.6){\alpha_1} &&&&c_1 \edge[dll]_{\gamma_1} \edge[drr]_(0.6){\gamma'_1} &&&& &&&a_0 \edge[dl]^{\alpha_0} &&&&& &&&a_0 \edge[dl] &&&&&&&&&  \\
d_0 &&&&a_0 \dashedge[dl]_{\alpha'_0} \dashedge[dr]^{\alpha_0} &&&&b_0 \dashedge[dl]_{\beta'_0} \dashedge[dr]^{\beta_0} &&&&c_0 \dashedge[dl]_{\gamma'_0} \dashedge[dr]^{\gamma_0} &&&&b_0 \dashedge[dl] \dashedge[dr] &&&&c_0 \dashedge[dl] \dashedge[dr] &&&&b_0 \dashedge[dl] \dashedge[dr] &&&&c_0 \dotedge[urr]  \dashedge[dl] \dotedge[drr]&&&&  \\
 &&&u &&c_0 \dashedge[dr]^{\gamma_0} &&v &&b_{-1} &&w &&c_{-1} &&v &&b_{-1} &&w &&c_{-1} &&v &&b_{-1} &&w &&&  \\
 &&& &&&c_{-1}
}$$
\centerline{\bf Orientation Diagram}

\bigskip

\noindent {\bf Claim 2.}  Let $m \ge2$.  Then every finite dimensional string module which belongs to  $\la_m\text{-}\mod$ but not to $\la_{m-1}\text{-}\mod$ has infinite projective dimension.  
\smallskip

\demo{Proof} Let $M = \St(w)$ be a string module in $(\la_m\text{-} \mod) \setminus (\la_{m-1}\text{-}\mod)$, where $w = p^{-1}_1 q_1 \dots p^{-1}_s q_s$ and the syllables $p_1^{-1}$ and $q_s$ may be trivial.  Once we have proved the claim for $m = 2$, it will follow for all $m > 2$ by induction:  Indeed, suppose $m \ge 3$.  In light of Proposition 2.1 and the graphs of the indecomposable projective $\la_m$-modules, the first syzygy $\Omega^1(M)$ is again a nonzero string module; it either belongs to $(\la_{m-1}\text{-} \mod) \setminus (\la_{m-2}\text{-}\mod)$, or else $m = 3$, and $\Omega_1(M)$ is isomorphic to one of the simple modules $S_{b_1}$, $S_{c_1}$.  (To verify this, observe that the only string modules $M$ in $(\la_{3}\text{-} \mod) \setminus (\la_{2}\text{-}\mod)$ which have a syzygy in $\la_1$-$\mod$ are of the form $\la x_3/  \soc \la x_3$ where $x \in \{a,b,c\}$.)

Next we focus on the case $m = 2$, where the argument is cumbersome, if elementary.  We include enough detail to indicate a strategy for the check: First, one verifies that $\pdim M = \infty$ whenever $M$ is a string module with simple top; in particular, this includes the simple modules $S_{a_2}$, $S_{b_2}$, $S_{c_2}$.  

Now assume that the top of $M$ is not simple. By hypothesis, $M$ has at least one top element normed by an idempotent among $e_{a_2}$,  $e_{b_2}$, $e_{c_2}$; in fact, either the first or the last nontrivial syllable of $w$ involves a path in $Q^{(2)}$ which does not belong to $Q^{(1)}$.  Since $\St(w) \cong \St(w^{-1})$, where $w^{-1} = q^{-1}_s p_s \dots q^{-1}_1 p_1$, it is therefore harmless to assume that either $p_1$ does not belong to $Q^{(1)}$ or else $p_1$ is trivial and $q_1$ does not belong to $Q^{(1)}$.

If $p_1 = \gamma'_2$, then $w = (\gamma'_2)^{-1} \gamma_2 (\beta'_2)^{-1} \cdots$, or $w = (\gamma'_2)^{-1} (\gamma_1 \gamma_2) (\alpha_1)^{-1} \cdots$.  In either case, the second syzygy of $M$ has a direct summand $S_{c_{-1}}$, clearly of infinite projective dimension.  If $p_1$ is trivial and $q_1 = \gamma'_2$, the $w = \gamma'_2 (\beta_2)^{-1} \cdots$, and $\Omega^2(M)$ has a uniserial direct summand with top $S_{b_0}$ and radical $S_v$, again of infinite projective dimension.  

If $p_1 = \gamma_2$, then $w = (\gamma_2)^{-1} \gamma'_2 (\beta_2)^{-1} \cdots$, which places a direct summand $S_u$ into $\Omega^2(M)$.  If $p_1$ is trivial and $q_1 = \gamma_2$, then $w = \gamma_2 (\beta'_2)^{-1} \cdots$, in which case the second syzygy of $M$ has a uniserial direct summand with top $S_{b_0}$ and radical $S_{b_{-1}}$. In either case, this shows the projective dimension of $\St(w)$ to be infinite. 

Suppose that $p_1 = \gamma_1 \gamma_2$.  Then $w = (\gamma_1 \gamma_2)^{-1} \gamma'_2 (\beta_2)^{-1} \cdots$, and the second syzygy of $M$ has a direct summand $S_w$. If, on the other hand $p_1$ is trivial and $q_1 =  \gamma_1 \gamma_2$, then $w = (\gamma_1 \gamma_2) (\alpha_1)^{-1} \cdots$, and we again find the uniserial module with top $S_{b_0}$ and radical $S_{b_{-1}}$ as a direct summand of $\Omega^2(M)$.

Similarly, one deals with the cases, where $\bullet$ $p_1 = \beta_2$ or  $\bullet$ $p_1$ is trivial and $q_1 = \beta_2$ or $\bullet$ $p_1 = \beta'_2$, or $\bullet$ $p_1$ is trivial and $q_1 = \beta'_2$.

The choice $p_1 = \gamma_2 \alpha'_2$ leads to a single word extending $p_1^{-1}$ to the right, namely $(\gamma_2 \alpha'_2)^{-1} \alpha_2$; this is due to the fact that $\alpha_2$ is the only arrow terminating in $a_1$.  The corresponding string module has simple top, and hence was addressed in the preliminary step.  If $p_1$ is trivial and $q_1 =  \gamma_2 \alpha'_2$, then $w$ is on the list  $\gamma_2 \alpha'_2 (\beta'_2)^{-1}$, $\gamma_2 \alpha'_2 (\beta'_2)^{-1} \beta_2$, $\gamma_2 \alpha'_2 (\beta'_2)^{-1} \beta_2 (\gamma'_2)^{-1}$, $\gamma_2 \alpha'_2 (\beta'_2)^{-1} \beta_2 (\gamma'_2)^{-1} \gamma_2$, $\dots$.  For each of these words, the inverse is among the words that have already been discussed.

For any word $w = (\alpha_2)^{-1} (\gamma_2 \alpha'_2) (\beta'_2)^{-1} \cdots$, the second syzygy $\Omega^2(\St(w))$ has a uniserial direct summand whose top is $S_{c_0}$ and whose radical is $S_{c_{-1}}$.  As for the case, where $p_1$ is trivial and $q_1 = \alpha_2$:  Observe that there is no multisyllabic word $w$ extending $(p_1)^{-1} q_1$, since $\alpha_2$ is the only arrow terminating in $a_1$.  

In case $p_1 = \alpha'_2$, the only word properly extending the syllable $p_1^{-1}$ to the right is $(\alpha'_2)^{-1} \alpha_2$, again leading to a string module with simple top. If, on the other hand $p_1$ is trivial and $q_1 = \alpha'_2$, the (essentially) monosyllabic word $p_1^{-1} q_1$ has no proper right extension, since the arrow $\alpha'_2$ is the only one terminating in $c_2$; thus it again results in a string module with simple top.  
This establishes Claim 2. \qed
\enddemo

\noindent {\bf Claim 3.} All band modules in $\la_{1}\text{-}\mod$ have projective dimension $1$.  For $m \ge 2$, all band modules in $(\la_{m}\text{-}\mod) \setminus (\la_{m-1}\text{-}\mod)$ have projective dimension $m$.  Combined with Theorem 0 and the preceding claims, this yields $\findim \la_m = r+1$ whenever $1 \le m \le r+1$.

\demo{Proof} We first derive the final statement from the first. That $\findim \la_m \ge \findim \la_1 = r+1$ is due to the fact that $\la_1$-$\mod$ is a subcategory of $\lamod$.  In light of Theorem 0, the reverse inequality will follow from Claim 2, once the postulated homological behavior of band modules is proved.   

Any band module $B$ in $\la_1$-$\mod$ is based on the primitive word $v_1 = (\beta_1)^{-1} \beta'_1 (\alpha'_0 \gamma_1)^{-1} \gamma'_1$ or, equivalently, on the words resulting from inversion or a permutation of the two pairs of syllables of $v_1$.  Hence, $\Omega^1(B)$ is a nonzero direct sum of copies of  $\la e_{b_0}$ and $\la e_{c_0}$, showing $B$ to have projective dimension $1$.  To deal with $m=2$, we note that all band modules over $\la_2$ have syzygies in $\la_1$-mod.  Indeed, the only primitive word (up to inversion and cyclic permutation of pairs of syllables) involving syllables not available in the alphabet for $\la_1$ is $v_2 = (\beta'_2)^{-1} \beta_2 (\gamma'_2)^{-1} \gamma_2$.  By Proposition 2.2,  the syzygies of the band modules based on this primitive word are again band modules; the orientation diagram shows these syzygies to be band modules in $\la_1$-mod in fact.  This implies $\pdim B = 2$, whenever $B$ is a band module in $\la_2\text{-} \mod \setminus \la_1\text{-}\mod$.

Finally suppose that $m \ge 3$.  Then the only primitive word $v_m$ not consisting of syllables defined over $\la_{m-1}$, again up to inversion and permutation of pairs of syllables, are as follows:  $v_3  = (\alpha'_2 \alpha_3)^{-1} \alpha'_3 (\gamma_3)^{-1} \gamma'_3$, and $v_m = (\alpha_m)^{-1} \alpha'_m (\beta_m)^{-1} \beta'_m$ for $m > 3$.  On inspection of the orientation diagram, we find that the (first) syzygy of any band module in $\la_m\text{-} \mod \setminus \la_{m-1}\text{-}\mod$ is a band module in $\la_{m-1}\text{-} \mod \setminus \la_{m-2}\text{-}\mod$, whence our assertion follows by induction. \qed 
\enddemo

We precede the final claim with an outline of the end game:  Clearly, the first syzygy of the module $\St(w)$, where $w = \alpha_1(\gamma_1)^{-1} \gamma'_1$, has a direct summand $S_{d_0}$ of projective dimension $r$.   In showing that, in contrast to Claim 3, the big finitistic dimensions of the $\la_m$ keep growing as $m$ increases, the finite words over $\la_1$ extending the word $w$ on the right  --  by way of the syllables $(\beta_1)^{-1}$, $\beta'_1$, $(\alpha_0 \gamma_1)^{-1}$, $\gamma'_1$  --  play the pivotal role.  From the orientation diagram one gleans that, as syllables are added, the first syzygies of the corresponding string modules retain the direct summand $S_{d_0}$, accrue additional projective summands in $\la_1$-$\mod$, and sport precisely one indecomposable direct summand of infinite projective dimension (rotating among finitely many isomorphism types).  This exceptional summand keeps moving to the right in the process, if the syzygies are displayed in the orientation of the diagram.  As a consequence, the interloper summand of infinite projective dimension disappears as we pass to the right infinite word that results from a cyclic repetition of the four add-on syllables.  The corresponding generalized string module is called $M_1$ in the proof of Claim 4.  The crucial point is that $M_1$ arises as an $(m-1)$-th syzygy of a module in $\la_m$-$\Mod$ for $m \ge 2$.  
\medskip

\noindent {\bf Claim 4.}   $\Findim \la_m = r + m$ for $1 \le m  \le r+1$.  In particular, $\Findim \la_{r+1} = 2r + 1$.

\demo{Proof}  We define $M_1 \in \la_1\text{-}\mod$ to be the generalized string module with graph
$$\xymatrixcolsep{0.5pc}\xymatrixrowsep{1.5pc}
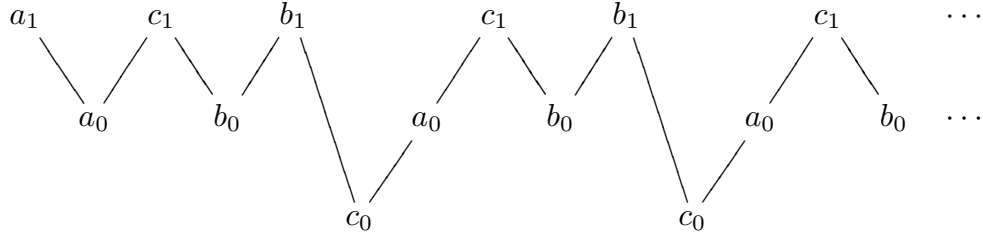
\xymatrix{
a_1 \edge[dr] &&c_1 \edge[dl] \edge[dr] &&b_1 \edge[dl] \edge[ddr] &&&c_1 \edge[dl] \edge[dr] &&b_1 \edge[dl] \edge[ddr] &&&c_1 \edge[dl] \edge[dr] &&\cdots  \\
 &a_0 &&b_0 &&&a_0 \edge[dl] &&b_0 &&&a_0 \edge[dl] &&b_0 &\cdots  \\
 &&&&&c_0 &&&&&c_0
 }$$
 The graph of this module can be found on the next-to-lowest level of the orientation diagram (single edges).

Next we define $M_2 \in \la_2\text{-}\mod$ as the generalized string module with graph
$$\xymatrixcolsep{0.5pc}\xymatrixrowsep{1.5pc}
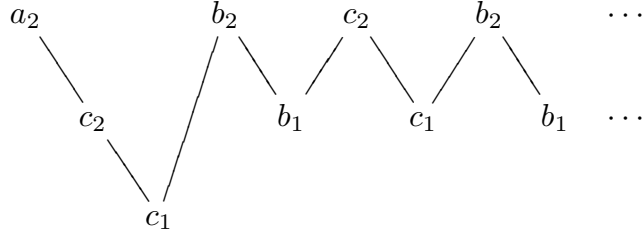
\xymatrix{
a_2 \edge[dr] &&&b_2 \edge[ddl] \edge[dr] &&c_2 \edge[dl] \edge[dr] &&b_2 \edge[dl] \edge[dr] &&\cdots  \\
 &c_2 \edge[dr] &&&b_1 &&c_1 &&b_1 &\cdots  \\
 &&c_1
 }$$
 It can be found on the level directly above $M_1$ in the orientation diagram (double edges).

$M_3 \in \la_3\text{-}\mod$ is the generalized string module with graph
$$\xymatrixcolsep{0.5pc}\xymatrixrowsep{1.5pc}
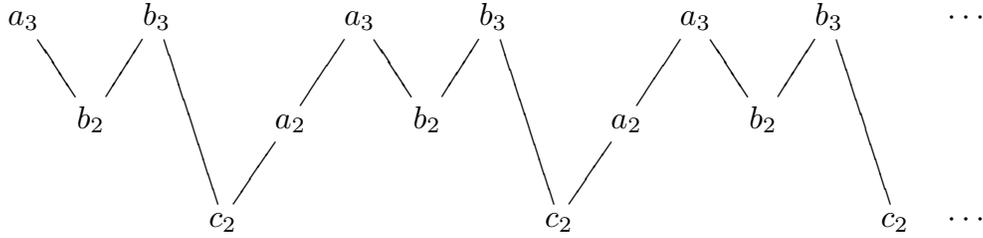
\xymatrix{
a_3 \edge[dr] &&b_3 \edge[dl] \edge[ddr] &&&a_3 \edge[dl] \edge[dr] &&b_3 \edge[dl] \edge[ddr] &&&a_3 \edge[dl] \edge[dr] &&b_3 \edge[dl] \edge[ddr] &&\cdots  \\
 &b_2 &&&a_2 \edge[dl] &&b_2 &&&a_2 \edge[dl] &&b_2  \\
 &&&c_2 &&&&&c_2 &&&&&c_2 &\cdots
 }$$
 Its graph can be retraced in the orientation diagram above $M_2$ (single edges).
 
 For $m \ge 4$, finally, $M_m \in \la_m\text{-}\mod$ is pinned down by the graph
$$\xymatrixcolsep{0.5pc}\xymatrixrowsep{1.5pc}
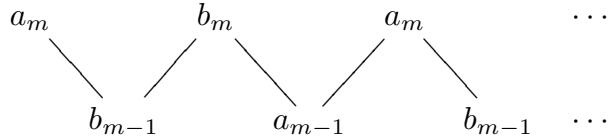
\xymatrix{
a_m \edge[dr] &&b_m \edge[dl] \edge[dr] &&a_m \edge[dl] \edge[dr] &&\cdots  \\
 &b_{m-1} &&a_{m-1} &&b_{m-1} &\cdots
 }$$
 \smallskip
 
To see that $\pdim M_1 = r + 1$, consult the diagram to obtain $\Omega^1(M_1) \cong S_{d_0} \oplus \la_2 e_{a_0} \oplus \la_2 e_{b_0}^{(\NN)} \oplus \la_2 e_{c_0}^{(\NN)}$.  Among the listed direct summands, the first has projective dimension $r$, while all others are projective.  The diagram further displays the fact that $\Omega^1(M_m) \cong M_{m-1}$ for all $m \ge 2$.  This yields $\pdim M_m = r + m$, and thus entails the inequality $\Findim \la_m \ge r+m$. 

For the reverse inequality, we first verify that $\Findim \la_2 \le r + 2$, keeping in mind that $\findim \la_1 = r+1$ by Claim 1. To that end, we show that any module $M \in \pinf(\la_2 \text{-}\Mod)$ has syzygy $\Omega^1(M) \in \la_1\text{-}\Mod$.  Indeed, if this were false, $M$ would have a top element normed by $e_{a_2}$ which is annihilated by $\alpha'_2$.  Since $\alpha_2$ is the only arrow in $Q^{(2)}$ ending in the vertex $a_1$, this would entail that $M$ has a direct summand which is either isomorphic to the simple module $S_{a_2}$ or to the uniserial  $\la_2$-module $\vcenter{ \xymatrixrowsep{2pc} \xymatrix{ a_2 \edge[d]\\ a_1 }}$.  In light of Claim 2, both of these possibilities are ruled out by finiteness of $\pdim M$.  Thus $\pdim \Omega^1(M) \le r+1$, and consequently $\pdim M \le r+2$.  

For every choice of $m \ge 3$, any $\la_m$-module has a syzygy in $\la_{m-1}\text{-}\Mod$, whence $M \in \pinf(\la_m \text{-}\Mod)$ entails $\pdim M \le r+ m$ by induction. This fills in the missing inequality $\Findim \la_m \le r+m$. \qed
\enddemo

Taking $\la$ to be $\la_{r+1}$, we thus obtain $\findim \la = r+1$ by Claim 3 and $\Findim \la = 2r + 1$ by Claim 4.

\proclaim{Theorem 3.1}  For any positive integer $r$, there exists a special biserial algebra $\la$ with the property that $\findim \la = r+1$, while $\Findim \la =  2r + 1 = 2 \findim \la - 1 $. \endproclaim

\definition{Remark}  For $m > r+1$, the little finitistic dimension of $\la_m$ starts increasing beyond that of $\la_1$.  Indeed, for all positive integers $m$, we have $\findim \la_m  \ge m$ by Claim 3.
\enddefinition

\enddocument